\documentclass[11pt,reqno,final]{amsart}
\newtheorem{theorem}{Theorem}[section]

\newtheorem{conjecture}[theorem]{Conjecture}

\newtheorem*{definition}{Definition}

\newtheorem{lemma}[theorem]{Lemma}

\newtheorem{proposition}[theorem]{Proposition}

\usepackage{url,upref,amscd}
\begin{document}
\begin{abstract}
 In this paper we study a key exchange protocol similar to the
 Diffie-Hellman key 
exchange protocol, using abelian subgroups of the automorphism group of a
non-abelian nilpotent group. We also generalize group no.92 of the Hall-Senior
table \cite{halltable} to an arbitrary prime $p$ and show that, for
those groups, 
the group of central automorphisms is commutative. We use these for the
key exchange we are studying.  
\end{abstract}
\title{The Diffie-Hellman key exchange protocol and non-abelian nilpotent groups}
\author{Ayan Mahalanobis}
\address{Department of Mathematical Sciences,
Stevens Institute of Technology,
Hoboken, NJ 07030, U.S.A.}
\email{amahalan@stevens.edu}
\today
\maketitle
\section{Introduction}
In this paper we generalize the Diffie-Hellman key exchange protocol from
a cyclic 
group to a finitely presented non-abelian nilpotent group of class $2$. 
Similar efforts were made in \cite{braid3,braid4,braid1} to use braid
groups,
a family of finitely presented non-commutative groups \cite{links,braid5}, in
key exchange. We also refer to \cite[Section 3]{sh1} for a formal
description of a key exchange protocol similar to
ours\begin{footnote}{The author expresses his gratitude to the referee
    for this reference.}\end{footnote}. 
Our efforts are not solely directed to construct an efficient and fast key
exchange protocol. We also try to understand the conjecture, \emph{the discrete
logarithm problem is equivalent to the Diffie-Hellman problem in a cyclic
group}. We develop and study protocols where, at least theoretically,
non-abelian groups can be used to share a secret or exchange private keys
between two people 
over an insecure channel. This development is significant
because nilpotent or, more specifically, $p$-groups have nice presentations and
computation in those groups is fast and easy \cite[Chapter 9]{sims}. So our
work can be seen as a nice 
application of the advanced and developed subject of $p$-groups and
computation with $p$-groups.

The frequently used public key cryptosystems are
slow and use mainly number theoretic complexity. The specific cryptographic
primitive that we 
have in mind is \emph{the discrete logarithm problem}, DLP for
short. DLP is general enough to be defined in an arbitrary cyclic group as
follows:  
let $G=\langle g\rangle$ be a cyclic group generated by $g$ and let
$g^n=h$, where $n\in\mathbb{N}$. Given 
$g$ and $h$, DLP is to find $n$ \cite[Chapter 6]{stinson}. The security of
the discrete logarithm 
problem depends on the representation of the group. It is trivial in
$\mathbb{Z}_n$, but is 
much harder (no polynomial time algorithm known)
in the multiplicative group of a finite field and even harder (no
sub-exponential time algorithm known) in the group of elliptic curves which
are not supersingular \cite{blakeCurve}. But with
the invention of sub-exponential algorithms for breaking the discrete logarithm
problem, like the index calculus and Coppersmith's algorithm, 
multiplicative groups of finite fields are no longer that attractive especially
the ones of characteristic $2$.

The discrete logarithm problem can be used in many other groups like
the group of
elliptic curves, in which case a cyclic group or a big enough cyclic 
component of an abelian group is used. In this article we propose a
generalization of DLP or more specifically the Diffie-Hellman key exchange
protocol in situations where the group has more than one
generator, i.e., in a finitely presented non-abelian group. Let $f$ be an
automorphism of a finitely presented group 
$G$ generated by $\{a_1,a_2,\ldots ,a_n\}$. If one knows the action of
$f$ on $a\in G$, i.e., $f(a)$, then it is difficult for him to tell the action
of $f$ on any other $b\in G$ i.e., $f(b)$. We describe this in detail later
under the name ``the general discrete logarithm problem''. In this paper we
work with finitely presented groups in terms of generators and relations and do
not consider any representation of that group. Though that seems
to be a good idea for future research.

Now suppose for a moment that $G=\langle g\rangle $ is a cyclic group and that
we are given $g$ and $g^n$ where $\gcd(n,|G|)=1$. DLP is to find $n$.
Notice that in this case the map $x\mapsto x^n$ is an automorphism. If we
conjecture that finding the automorphism is finding
$n$ then one way to see DLP, in terms of group theory, is to find
the automorphism from its image on one element. This is the central idea that
we want to generalize to non-abelian finitely presented groups,
especially to a family of $p$-groups of class $2$. This explains our
choice of the name the general discrete logarithm problem.  

To work with a finitely presented group and its automorphisms the
following properties of the group are needed.
\begin{itemize}
\item A consistent and natural representation of the elements in the group.
\item Computation in the group should be fast and easy.
\item The automorphism group should be known and the automorphisms should have
  a nice enough presentation so that images can be computed quickly.
\end{itemize}
We note at this point that for a $p$-group the first two requirements are
satisfied \cite[Chapter 9]{sims}. 
\section{Our Contribution in this article}
The central idea behind this article is to study a generalization of
the discrete logarithm problem (DLP) that we call the general discrete
logarithm problem (GDLP). As a cryptographic primitive the concept of
GDLP seems to be secure (see Section 4.1).

To use GDLP we use a Diffie-Hellman like key exchange protocol using
finitely presented $p$-groups with an abelian central automorphism
group. In this case the security depends not only on GDLP but also on
GDHP (see Section 4.2) which turns
out to be insecure in the specific case we are studying.

Section 8 of this paper contains a brief survey of all the group
theoretic results necessary for a reader to understand the later part
of this paper. However, a knowledgeable reader might choose to ignore
Section 8 altogether and come back to it when required. In Section 10
we survey the existing literature for groups with abelian
automorphism group and show that none of them are adequate for the key
exchange we are studying. 

We found no groups readily available in the literature, hence we had to
develop a family of groups $G_n(m,p)$ with abelian central
automorphism group (Section 10). This is a significant contribution to
the theory of finite groups because $G_n(m,p)$ is a generalization of group
no.~92 of the Hall-Senior table. We describe the group of
automorphisms for this group and further prove that this group is
Miller if and only if $p=2$.

We do not claim that the key exchange protocol is
secure. Rather, we show that the key exchange protocol is
insecure for the particular family of  groups that we picked. Our
study raises two important questions which are of interest both
mathematically as well as cryptographically. 
\begin{description}
\item[a] Are there groups different from $G_n(m,p)$, with an abelian central
  automorphism group, for which the key exchange protocol is secure?
\item[b] Does there exist any cryptographic protocol with reductionist
  security proof, where the security of the protocol depends only on the
  discrete logarithm problem? If one can find such a protocol using
  cyclic groups then that could be generalized using GDLP, and since
  we claim that GDLP is a secure primitive, this will give rise to a
  secure cryptosystem using non-abelian groups.
\end{description}  
\section{Some notations and Definitions}
We now describe some of the definitions and notations that will be used in this
paper. The notations used are standard:
\begin{itemize}
\item $G$ will denote a finite group. $Z=Z(G)$ denotes the center of the group
  $G$ and will be denoted by $Z$ if no confusion can arise.
\item $G^\prime=[G,G]$ is the commutator subgroup of $G$.
\item $\text{Aut}(G)$ and $\text{Aut}_c(G)$ are the group of automorphisms and
  the group of central automorphisms of $G$, respectively. 
\item $\Phi(G)$ is the Frattini subgroup of $G$, which is the intersection of
  all maximal subgroups of $G$.
\item We denote the commutator of $a,b$ by $[a,b]$ where
  $[a,b]=a^{-1}b^{-1}ab$. 
\item The exponent of a $p$-group $G$, denoted by exp$(G)$, is the largest
  power of $p$ that is the order of an element in $G$.
\end{itemize}
The following commutator formulas hold for any element $a$,$b$ and $c$ in any
group $G$.
\begin{description}
\item[(a)] $a^b=a[a,b]$
\item[(b)] $[ab,c]=[a,c]^b[b,c]=[a,c][a,c,b][b,c]$ it follows that in a
  nilpotent group of class $2$, $[ab,c]=[a,c][b,c]$
\item[(c)] $[a,bc]=[a,c][a,b]^c=[a,c][a,b][a,b,c]$ it follows that in a
  nilpotent group of class $2$, $[a,bc]=[a,b][a,c]$
\item[(d)] $[a,b]^{-1}=[b,a]$  
\end{description}
The proof of these formulas follow from direct computation or can be found in
\cite{khukhro}.
\begin{definition}[Miller Group]
A group $G$ is called a Miller group if it has an abelian automorphism
group, in 
other words, if $\text{Aut}(G)$ is commutative then the group $G$ is Miller.
\end{definition}
\begin{definition}[Central Automorphisms]
Let $G$ be a group, then $\phi\in\text{Aut}(G)$ is called a central
automorphism if $g^{-1}\phi(g)\in Z(G)$ for all $g\in G$. Alternately, one
might say that 
$\phi$ is a central automorphism if $\phi(g)=gz_{\phi,g}$ where
$z_{\phi,g}\in Z(G)$ depends on $g$ and $\phi$. If $\phi$ is clear from the
context then we can simplify the notation as $\phi(g)=gz_g$. 
\end{definition}
Apart from inner automorphisms, central automorphisms are second best in terms
of nice description. They are very attractive for cryptographic
purposes, since it is easy to describe the automorphisms and compute the image
of an arbitrary element.
\begin{theorem}
The centralizer of the group of inner automorphisms is the group of central
automorphisms. Moreover a central automorphism fixes the commutator
elementwise.    
\end{theorem}
This theorem first appears in \cite{fournelle} which refers to \cite{hall} and
\cite{zeanhaus}.

\begin{definition}[Polycyclic Group]
Let $G$ be a group, a finite series of subgroups in $G$
\[G=G_0\trianglerighteq G_1\trianglerighteq G_2\trianglerighteq
G_3\trianglerighteq\ldots\trianglerighteq G_n=1\]
is a polycyclic series if $G_i/G_{i+1}$ is cyclic and $G_{i+1}$ is a normal
subgroup of $G_i$. Any group with polycyclic series is a polycyclic group. 
\end{definition}
It is easy to prove that finitely generated nilpotent
groups are polycyclic, hence any finitely generated $p$-group is
polycyclic. Let $a_i$ be an 
element in $G_i$ whose image generates $G_i/G_{i+1}$. Then the sequence
$\{a_1,a_2,\ldots ,a_n\}$ is called a polycyclic generating set. It is easy to
see that $g\in G$ can be written as $g=a_1^{\alpha_1}a_2^{\alpha_2}\ldots
a_n^{\alpha_n}$, where $\alpha_i$ are integers. If
$g=a_1^{\alpha_1}a_2^{\alpha_2}\ldots a_n^{\alpha_n}$ where $0\leq\alpha_i<
m_i$, $m_i=\left|G_i:G_{i+1}\right|$ then the expression is a collected
word. Each element $g\in G$ can be expressed by a unique collected
word. Computation with these collected words is easy and implementable in
computer, for more information on this topic see \cite[Section
9.4]{sims} and also \cite[polycyclic package]{GAP4}.  
\section{Key Exchange}
We want to follow the Diffie-Hellman Key exchange protocol using a commutative
subgroup of the automorphism group of a finitely presented group $G$. The
security of the Diffie-Hellman key exchange protocol in a cyclic group rests on the following three
factors: 
\begin{description}
\item[DLP] The discrete logarithm problem.
\item[DHP] The Diffie-Hellman problem.
\item[DDH] The decision Diffie-Hellman problem
  \cite{blake1,boneh,galb1,igor,gonzalez}. 
\end{description}
We have already described the discrete logarithm problem. The Diffie-Hellman
problem is the following:
let $G=\langle g\rangle$ be a cyclic group of order $n$. One knows $g$, $g^a$
and $g^b$, and the problem is to compute $g^{ab}$. It is not known if DLP is
equivalent to DHP. The decision Diffie-Hellman problem is more
subtle. Suppose that DHP is a hard problem, so it is impossible to compute
$g^{ab}$ from $g^a$, $g^b$ and
$g$. But what happens if someone can
compute or predict $80\%$ of the
binary bits of $g^{ab}$ from
$g^a$, $g^b$ and $g$, then the adversary will have 
 $80\%$ of the shared secret or the
 private key; that is most of the private key. This is clearly
 unacceptable. It is often hard to formalize DDH in exact mathematical
 terms (\cite[Section 3]{boneh}); the best formalism offered is a randomness
 criterion for the bits of the key. In DDH we ask the question, given the
 triple $g^a, g^b$ and $g^c$ is $c=ab\mod n$?
 But there is no known link between DDH and any mathematically hard problem
 for the Diffie-Hellman key exchange protocol in cyclic groups. 

Clearly, solving the
 discrete logarithm problem solves the Diffie-Hellman problem and solving
 the Diffie-Hellman problem solves the decision Diffie-Hellman problem.

As is usual, we denote by Alice and Bob, two people trying
to set up a private key over an insecure channel to communicate securely and
Oscar an eavesdropping adversary. In this paper the shared secret or the
private key is an element of a finitely presented group $G$.
\subsection{General Discrete Logarithm Problem} Let $G=\langle
a_1,a_2,\ldots,a_n\rangle$ and $f:G\rightarrow G$ be a
non-identity automorphism. Suppose one knows $f(a)$ 
and $a\in G$  then GDLP is to find $f(b)$ for any $b$ in $G$. 
Assuming the word problem is easy or presentation of the group is by means of
generators, GDLP is equivalent to finding $f(a_i)$ for all $i$ which in terms
gives us a complete knowledge of the automorphism. So in other words the
cryptographic primitive GDLP is equivalent to, ``\textit{finding the
  automorphism $f$ from the action of $f$ on only one element}''. 
\subsection{General Diffie-Hellman Problem} Let $\phi,\psi:G\rightarrow G$ be
arbitrary automorphisms such that $\phi\psi=\psi\phi$, and assume one
knows $a$, $\phi(a)$ and 
$\psi(a)$. Then GDHP is to find $\phi(\psi(a))$. Notice that GDHP is a
restricted form of GDLP, because in case of GDHP one has to compute
$\phi(\psi(a))$ for some fixed $a$, not $\phi(b)$ for an arbitrary $b$ in $G$. 
 There is an interesting GDHP attack due to Vladimir
 Shpilrain. To mount this attack one need not find $\phi$ but finds
 another automorphism $\phi^\prime$ such that
 $\phi^\prime\psi=\psi\phi^\prime$ and $\phi^\prime(a)=\phi(a)$. Since
 $\phi(\psi(a))=\psi(\phi^\prime(a))=\phi^\prime(\psi(a))$, the
 knowledge of the $\phi^\prime$ breaks the system. We will refer to
 this attack as the Shpilrain's attack. 

We now describe two key exchange protocols and do some
cryptanalysis. We denote by $G$ a finitely presented group and $S$ an abelian
subgroup of $\text{Aut}(G)$. 
\section{Key Exchange Protocol I}\label{keyexchange1}
Alice and Bob want to set up a private key. They select a group $G$
and an element $a\in G\setminus Z(G)$ over an insecure channel. Then Alice
picks a random automorphism $\phi_A\in S$ and sends Bob $\phi_A(a)$. Bob
similarly picks a random automorphism 
$\phi_B\in S$ and sends Alice $\phi_B(a)$. Both of them can now compute
$\phi_A(\phi_B(a))=\phi_B(\phi_A(a))$ which is their private key for a
symmetric transmission.
\begin{description}
\item[Step 1] Alice and Bob selects the group $G$ and an element $a\in
  G\setminus Z(G)$ in public. Notice that $G$ and $a$ are public
  information.
\item[Step 2] Alice and Bob picks,
  at random, two automorphisms $\phi_A$
  and $\phi_B$ from $S$ respectively. Notice that $\phi_A$ and
  $\phi_B$ are private information.
\item[Step 3] Alice and Bob compute $\phi_A(a)$ and $\phi_B(a)$
  respectively and exchanges them. Notice that $\phi_A(a)$ and
  $\phi_B(a)$ are public information.
\item[Step 4] Both of them compute
  $\phi_A\left(\phi_B(a)\right)=\phi_B\left(\phi_A(a)\right)$ from their private information; which is their private key.
\end{description}
\subsection{Comments on Key Exchange Protocol I}
Though initially it might seem that we do not have enough information to know
the automorphisms $\phi_A$ and $\phi_B$, it turns out that if we are
using automorphisms which fix conjugacy classes, like inner automorphisms,
then the security of the above scheme actually rests on the conjugacy problem.

Let $\phi_A(a)=x^{-1}ax$ and $\phi_B(a)=y^{-1}ay$ for some $x$ and $y$. Then
$\phi_A(\phi_B(a))\\=(yx)^{-1}a(yx)$. Since $a$, $\phi_A(a)$ and
$\phi_B(a)$ are known, if the conjugacy problem is easy in the group
then anyone can find $x$ and $y$ and break the system.

In the above scheme Oscar knows $G$ and $a$. If the automorphisms are central
automorphisms, then
he also sees $\phi_A(a)=az_{\phi_A ,a}$ and $\phi_B(a)=az_{\phi_B ,a}$. Oscar
can compute $z_{\phi_A, a}$ and
$z_{\phi_B ,a}$. Now if $G$ is a
special $p$-group 
($G^\prime=Z(G)=\Phi(G)$) then $Z(G)$ is fixed elementwise by both $\phi_A$
and $\phi_B$. Then
\begin{equation}\label{center attack}
\phi_A(\phi_B(a))=\phi_A(az_{\phi_B ,a})=az_{\phi_A ,a}z_{\phi_B ,a}.
\end{equation}
Oscar knows $a$ and can compute $z_{\phi_A ,a}$ and $z_{\phi_B ,a}$ and can
find the 
private key $\phi_A(\phi_B(a))$. In the literature all examples of Miller 
$p$-group with odd prime $p$ are special and the above key exchange is
fatally flawed for those groups.
\section{Key Exchange Protocol II}\label{keyexchange2}
In this case Alice and Bob want to set up a private key and they set up a
group $G$ over an insecure channel. Alice chooses a random non-central
element $g$ and a random automorphism $\phi_A\in S$ and sends Bob
$\phi_A(g)$. Bob 
picks another automorphism $\phi_B\in S$ and computes $\phi_B(\phi_A(g))$ and
sends it back to Alice. Alice, knowing $\phi_A$, computes $\phi_A^{-1}$ which
gives her $\phi_B(g)$ and picks
another random automorphism $\phi_H\in S$ and computes 
$\phi_H(\phi_B(g))$ and sends it back to Bob. Bob, knowing $\phi_B$ computes
$\phi^{-1}_B$ which gives him
$\phi_H(g)$ which is their private key. Notice that Alice never
reveals $g$ in public.

\begin{description}
\item[Step 1] Alice and Bob set up
  the group $G$. Notice that $G$ is
  public information.
\item[Step 2] Alice picks $g\in
  G\setminus Z(G)$ and a random $\phi_A\in
  S$. Then she computes $\phi_A(g)$ and sends that to Bob. Notice that
  $g$ and $\phi_A$ are private but $\phi_A(g)$ is public.
\item[Step 3] Bob picks $\phi_B\in S$ at random and computes
  $\phi_B\left(\phi_A(g)\right)$ and sends that back to Alice. Notice
  that $\phi_B$ is private but $\phi_B\left(\phi_A(g)\right)$ is
  public.
\item[Step 4] Alice computes $\phi_A^{-1}$ and then computing
  $\phi_A^{-1}\left(\phi_B\left(\phi_A(g)\right)\right)$ she gets
  $\phi_B(g)$. 
\item[Step 5] Alice now picks another random automorphism $\phi_H\in
  S$ and computes $\phi_H\left(\phi_B(g)\right)$ and $\phi_H(g)$. She
  then sends $\phi_H\left(\phi_B(g)\right)$ to Bob but keeps
  $\phi_H(g)$ private.
\item[Step 6] Similar to Step 4, Bob computes $\phi_H(g)$. Now both
  Alice and Bob knows $\phi_H(g)$ and it is their common key.
\end{description}  
\subsection{Comments on Key Exchange Protocol II}
Notice that for central automorphisms, $\phi_A$ and $\phi_B$, $\phi_A(g)=gz_{\phi_A ,g}$; since $g$ is not known Oscar doesn't know
$z_{\phi_A ,g}$ but if $G$ is special ($Z(G)=G^\prime=\Phi(G)$) then
$\phi_B(gz_{\phi_A ,g})\\=gz_{\phi_B ,g}z_{\phi_A ,g}$ from which
$z_{\phi_B ,g}$ can be computed. Now
$\phi_H(\phi_B(g))=gz_{\phi_B,g}z_{\phi_H ,g}$ is a public
information; so using $z_{\phi_B ,g}$ one 
can compute $gz_{\phi_H ,g}$, which is $\phi_H(g)$ and the scheme is
broken. As one 
clearly sees, this attack is not possible if the group is not special.

The reader might have noticed at this point that all the attacks are
GDHP. So certainly in some groups GDHP is easy, even though GDLP is hard.

As we know, any automorphism in $G$ can be
seen as a restriction of an inner automorphism in $\text{Hol}(G)$
(see \cite{kurosh,maria} for further details on the holomorph of a
group). Solving the  
conjugacy problem in $\text{Hol}(G)$ will break the key exchange
protocols for any 
automorphism. On the other hand, operation in $\text{Hol}(G)$ is twisted so it
is possible that the conjugacy problem in $\text{Hol}(G)$ is difficult even
though it is easy in $G$. Since any cyclic group is a Miller group,
success of the 
holomorph attack would prove insecurity in DLP. Therefore we believe that the
holomorph attack will not be successful in many cases. Though more work needs
to be done on this.
\section{Key Exchange using Braid Groups}
In \cite{braid1} a similar key exchange protocol was defined,
in this section we mention some similarities of their approach to ours.
We also mention how our system
generalizes their system which uses braid
      groups. See also \cite{polybraid}.

We define braid group as a finitely presented group, though there are
fancy pictorial ways to look at braids and multiplication of braids. An
interested reader can look in
\cite{links,braid5}. The braid
group $B_n$ with $n$-strands is defined as:
\begin{eqnarray*}
B_n=\left\langle
  \sigma_1,\ldots,\sigma_{n-1}:\sigma_i\sigma_j\sigma_i=\sigma_j
\sigma_i\sigma_j\,\text{if}\,|i-j|=1, 
\sigma_i\sigma_j=\sigma_j\sigma_i \,\text{if}\, |i-j|\geq 2
  \right\rangle 
\end{eqnarray*}
In \cite{braid1}, the authors found two subgroups $A$ and $B$ of the
group of inner automorphisms of $B_n$, $\text{Inn}(B_n)$, such that,
if $\phi\in A$ and $\psi\in B$, then
$\phi(\psi(g))=\psi(\phi(g))$ for
$g\in B_n$. Then the key exchange proceeds
similar 
to the Key Exchange Protocol I above; with the restriction that Alice chooses
automorphisms from A and Bob chooses automorphisms from B. There is also a
different approach to key exchange using braid groups as in \cite{braid3,braid4}.

In the same spirit as \cite{braid1} we can develop a 
key exchange protocol similar to the key exchange protocol I, where we take two
subgroups $A$ 
and $B$ in $\text{Aut}(G)$ such that for $\phi\in A$ and $\psi\in B$,
$\phi(\psi(g))=\psi(\phi(g))$ for all $g\in G$. The use of inner automorphisms is only possible
when the conjugacy or the generalized conjugacy problem (conjugator search
problem) is known to be hard.

There are significant differences in our approach to that of the approach
in \cite{braid1}. In \cite{braid1}, the authors choose a group and then try to
use that group in cryptography. On
the other hand, we take the fundamental
concept as the discrete logarithm problem, generalize it using
automorphisms of a non-abelian group and then look for groups favorable to
us. The fact that the central idea in braid group key exchange turns out to be
similar to ours is encouraging.

It is intuitively clear at this point that we should start looking
for groups with abelian automorphism group, i.e., Miller groups.
\section{Some useful facts from group theory}
The term Miller Group is not that common in the literature. It was 
introduced by Earnley in \cite{earnley}. Miller was the first to study groups
with abelian automorphism group in \cite{miller}.
Cyclic groups are good examples of Miller groups. G.A. Miller also proved that
no non-cyclic abelian group is Miller. 

Charles Hopkins began a list of necessary conditions for a Miller group in
1927 \cite{hopkins}. He complained that very little is known about
those groups. The same is true today. Except for some sporadic examples of
groups with abelian automorphism groups, there is no sufficient condition
known for a group to be Miller. 

We now state some known facts about Miller groups which are
 available in the literature and which we shall need later. For proof of these
 theorems which we present in 
 a rapid fire fashion, the reader can look in any standard text books, like
 \cite{khukhro,rotman}, or the references there. 
\begin{proposition}
If $G$ is a non-abelian Miller group, then $G$ is nilpotent and of class $2$.
\end{proposition}
\begin{proof}
It follows from the fact that the group of inner automorphisms commute and
$G/Z(G)\cong\text{Inn}(G)$. 
\end{proof}
 Since a nilpotent group is a direct product of its Sylow $p$-subgroups $S_p$,
 and $\text{Aut}(A\times B)=\text{Aut}(A)\times\text{Aut}(B)$ whenever $A$ and
 $B$ are of relatively prime order, it is enough to study
 Miller $p$-groups for prime $p$.
\begin{proposition}
If $G$ is a $p$-group of class $2$, then
$\text{exp}(G^\prime)=\text{exp}(G/Z(G))$. 
\end{proposition}   
\begin{proposition}
In a $p$-group of class $2$,
$(xy)^n=x^ny^n[y,x]^{\frac{n(n-1)}{2}}$. Furthermore if
$\text{exp}(G^\prime)=n$ is odd, then $(xy)^n=x^ny^n$.
\end{proposition}
By definition, in a Miller group all automorphisms commute. Since central
automorphisms are the centralizer of the group of inner automorphisms, we have
proved the following theorem.
\begin{theorem}
In a Miller group $G$, all automorphisms are central.
\end{theorem}
It follows that to show a group is not Miller, all we have to do is
to produce a non-central automorphism.
\begin{proposition}
If the commutator and the center coincide then every pair of central
automorphisms commute.
\end{proposition}
\begin{proof}
Let $G$ be a group such that $G^\prime=Z(G)$. Then let $\phi$ and $\psi$ be
central automorphisms given by $\phi(x)=xz_{\phi ,x}$ and $\psi(x)=xz_{\psi
  ,x}$ where $z_{\phi ,x},z_{\psi ,x}\in G^\prime$. Then
$$\psi(\phi(x))=\psi(xz_{\phi ,x})=\psi(x)z_{\phi ,x}=xz_{\psi ,x}z_{\phi
  ,x}=xz_{\phi ,x}z_{\psi ,x}=\phi(\psi(x)).$$ 
\end{proof}
\begin{definition}[Purely non-abelian group]
A group $G$ is said to be a purely non-abelian group (PN group for short) if
whenever $G=A\times B$ where $A$ and
$B$ are subgroups of $G$ with $A$ abelian, then $A=1$. Equivalently $G$ has
no non-trivial abelian direct factor.  
\end{definition}
Let $\sigma:G\rightarrow G$ be a central automorphism. Then we define a map
$f_{\sigma}:G\rightarrow Z(G)$ as follows:
$f_{\sigma}(g)=g^{-1}\sigma(g)$. 
Clearly this map defines a homomorphism. 
The map $\sigma\mapsto f_{\sigma}$ is clearly a one-one map.
Conversely, if $f\in\text{Hom}(G,Z(G))$ then we define a map
$\sigma_f(g)=gf(g)$, $x\in G$. 
Clearly $\sigma_f$ is an endomorphism. It is easy to see that 
\[\text{Ker}(\sigma_f)=\{x\in G:\;f(x)=x^{-1}\}.\]
Hence it follows that $\sigma_f$ is an automorphism if and only if $f(x)\neq
x^{-1}$ for all $x\in G$ with $x\neq 1$.
\begin{theorem}\label{correspondence}
In a purely non-abelian group $G$, the correspondence $\sigma\rightarrow
f_{\sigma}$ is a one-one map of $\text{Aut}_c(G)$ onto $\text{Hom}(G,Z(G))$
\end{theorem}
\begin{proof}
See \cite{adney}.    
\end{proof}
For any $f\in\text{Hom}(G,Z(G))$ there is a map
$f^\prime\in\text{Hom}(G/G^\prime,Z(G))$ since $f(G^\prime)=1$. Furthermore,
corresponding to $f^\prime\in\text{Hom}(G/G^\prime,Z(G))$ there is
a map $f:G\rightarrow Z(G)$ explained in the following diagram
\[\begin{CD}
G @>\eta>> G/G^\prime @>f^\prime>> Z(G)
\end{CD}\]
where $\eta$ is the natural epimorphism.

Let $G$ be a $p$-group of class $2$, such that exp$(Z(G))=a$,
exp$(G^\prime)=b$ and exp$(G/G^\prime)=c$ and let $d=\min(a,c)$. Now
from the fundamental theorem of abelian groups, let
$$G/G^\prime =A_1\oplus A_2\oplus\ldots A_r\;\; \text{where}\;\; A_i=\langle
a_i\rangle$$
$$Z(G)=B_1\oplus B_2\oplus\ldots B_s\;\; \text{where}\;\; B_i=\langle
b_i\rangle$$
$r,s\in\mathbb{N}$ be the direct decomposition of $G/G^\prime$ and $Z(G)$. 
If the cyclic component $A_k=\langle a_k\rangle$ has exponent greater or
equal to the exponent of $B_j=\langle b_j\rangle$, then one can define a
homomorphisms $f:G/G^\prime\rightarrow Z(G)$ as follows
\[f(a_i)=\left\{
\begin{array}{lll}
b_j\;\;\text{where}\;\;i=k\\
1 \;\;\;\text{where}\;\;i\neq k
\end{array}
\right.
\]
From this discussion it is clear that for $f\in\text{Hom}(G,Z(G))$, $f(G)$
generates the subgroup 
\[\mathcal{R}=\{z\in Z(G):\; |z|\leq p^d,\; d=\min(a,c)\}.\]
\begin{definition}[Height]
In any abelian $p$-group $A$ written additively, there is a descending
sequence of subgroups 
\[A\supset pA\supset p^2A\supset \ldots \supset p^nA\supset
p^{n+1}A\supset\ldots\]
Then $x\in A$ is of height $n$ if $x\in p^nA$ but not in $p^{n+1}A$. In other
words the elements of height $n$ are those that drop out of the chain in the
$(n+1)^\text{th}$ inclusion.
\end{definition}
For further information on height see \cite{kaplansky}.

Since for a class $2$ group we have
\[\text{exp}(G/G^\prime)\geq\text{exp}(G/Z(G))=\text{exp}(G^\prime)\]   
it follows that $c\geq b$. Hence if $d=\min(a,c)$ then either $d=b$ or $d>b$. 
 
Let height$(xG^\prime)\geq b$, then $xG^\prime=y^{p^b}G^\prime$ for some
$y\in G$. Then for any $F\in\text{Hom}(G,G^\prime)$, $F(yG^\prime)^{p^b}=1$
implying $xG^\prime\in F^{-1}(1)$.
Conversely, let height$(xG^\prime)<b$. Then from the previous discussion it
is clear that there is a $F^\prime\in\text{Hom}(G/G^\prime,G^\prime)$ such that
$xG^\prime$ is not in the kernel, consequently there is a
$F\in\text{Hom}(G,G^\prime)$ such that $x\notin\text{ker}(F)$.
Combining these two facts we see that: 
\[\mathcal{K}=\bigcap\limits_{F\in\text{Hom}(G,G^\prime)}F^{-1}(1)=\left\{x\in 
G:\;\text{height}(xG^\prime)\geq b\right\}\]
\begin{proposition}
$\mathcal{K}\subseteq\mathcal{R}$
\end{proposition}
\begin{proof}
In a class $2$ group, if $x\in\mathcal{K}$ then $xG^\prime=y^{p^b}G^\prime$
for some $y\in G$ and exp$(G/Z)=b$ and $G^\prime\subseteq Z(G)$, hence $x\in
Z(G)$. 

Let $x\in\mathcal{K}$, then height$(xG^\prime)\geq b$, hence there is a
$y\in G$ such that $y^{p^b}G^\prime=xG^\prime$ i.e., $x=y^{p^b}z$ where
$z\in G^\prime$ and $y^{p^c}\in G^\prime$ and $c\geq b$. We have 
\[x^{p^c}=(y^{p^b})^{p^c}z^{p^c}=(y^{p^c})^{p^b}=1\]
Hence $|x|\leq\min(p^a,p^c)$ which implies that $x\in\mathcal{R}$.  
\end{proof}
\begin{proposition}
For a PN group $G$ of class $2$, if $\text{Aut}_c(G)$ is abelian then
$\mathcal{R}\subseteq\mathcal{K}$. 
\end{proposition}
\begin{proof}
In a PN group, using Theorem \ref{correspondence} and the notation there, two
central automorphisms $\sigma$ and $\tau$ 
commute if and only if $f_{\sigma},f_{\tau}\in\text{Hom}(G,Z(G))$
commute. Then for any $f\in\text{Hom}(G,Z(G))$ and
$F\in\text{Hom}(G,G^\prime)$ we have that $f\circ 
F=F\circ f=1$. Since $f(G^\prime)=1$, clearly $F\circ f(G)=1$ proving that
$\mathcal{R}\subseteq\mathcal{K}$.  
\end{proof}
Combining the above two propositions, we just proved that in a PN group
$G$ of class $2$, if
$\text{Aut}_c(G)$ is abelian then $\mathcal{R}=\mathcal{K}$. As discussed
earlier there are two cases $d=b$ and $d>b$. Adney and Yen proves that: 
\begin{proposition}
If $G$ is a non-abelian $p$ group of class $2$, and $\text{Aut}_c(G)$ is
abelian with $d>b$, then $\mathcal{R}/G^\prime$ is cyclic.
\end{proposition}
\begin{proof}
See \cite[Theorem 3]{adney}.
\end{proof}
\begin{theorem}[Adney and Yen]\label{adney}
Let $G$ be a purely non-abelian group of class $2$, $p$ odd, let
$G/G^\prime=\prod\limits_{i=1}^n\{x_iG^\prime\}$. Then the group
$\text{Aut}_c(G)$ is abelian if and only if 
\begin{itemize}
\item[(i)] $\mathcal{R}=\mathcal{K}$
\item[(ii)] either $d=b$ or $d>b$ and $\mathcal{R}/G^\prime=\{x_1^{p^b}G^\prime\}$
\end{itemize}
\end{theorem}
\begin{proof}
See \cite[Theorem 4]{adney}.
\end{proof}
From the proof of Proposition 8.5 it follows that in a group $G$ with
$Z(G)\leq G^\prime$, the central automorphisms commute.
\begin{theorem}\label{sanders}
The group of central automorphisms of a $p$-group $G$, where $p$ is odd, is a
$p$-group if and only if $G$ has no non-trivial abelian direct factor.
\end{theorem}
\begin{proof}
See \cite[Theorem B]{sanders} and its corollary.
\end{proof}
At this point we concentrate on building a cryptosystem. We note that Miller
groups in particular have no advantage over groups with abelian
central automorphism group. It is hard to construct Miller groups
and there is no known Miller group
for an odd prime, which is not special. So we
now turn towards a group $G$ such that $\text{Aut}(G)$ is not abelian but
$\text{Aut}_c(G)$ is abelian. We propose to use $\text{Aut}_c(G)$ rather than
$\text{Aut}(G)$ in the key exchange protocols described earlier.
\section{Signature Scheme based on conjugacy problem}
Assume that we are working with a group $G$ with commuting inner
automorphisms.

Alice publishes $\alpha$ and $\beta$ where $\beta=a^{-1}\alpha a$ and
keeps $a$ a secret. To sign a text $x\in G$ she picks an arbitrary
element $k\in G$ and computes $\gamma=k\alpha k^{-1}$ and then computes
$\delta$ such that $x=(\delta k)(a\gamma)^{-1}$. Now notice that 
\begin{eqnarray*}
x\alpha x^{-1}=&(\delta k)(a\gamma)^{-1}\alpha((\delta
k)(a\gamma)^{-1})^{-1}&\\
=&(\delta k)\gamma^{-1}a^{-1}\alpha a\gamma k^{-1}\delta^{-1}&\\ 
=&\delta\gamma^{-1}a^{-1}k\alpha k^{-1}a\gamma\delta^{-1}&\text{Inner
  automorphisms commute}\\
=&\delta\gamma^{-1}a^{-1}\gamma a\gamma\delta^{-1}\\
=&\delta a^{-1}\gamma a\delta^{-1}&\\
=&\delta (k\beta k^{-1})\delta^{-1}& \gamma=k\alpha k^{-1}\Rightarrow
a^{-1}\gamma a=k\beta k^{-1}\\
\end{eqnarray*}  
So to sign a message $x\in G$ Alice computes $\delta$ as mentioned and sends
$x,(k\delta)$. To verify the message one computes $L=x\alpha x^{-1}$ and
$R=\delta k\beta(\delta k)^{-1}$. If $L=R$ then the message is authentic
otherwise not.

There is a similar signature scheme in \cite{braid2}, where they exploit the
gap between the computational version (conjugacy problem) and the decision
version of the conjugacy problem (conjugator search problem) in braid
groups. We followed the El-Gamal signature scheme closely 
\cite[Chapter 7]{stinson}.
\subsection{Comments on the above Signature Scheme}
If one can solve conjugacy problem in the group then from the public
information $\alpha$ and $\beta$ he can find out $a$ and our scheme is
broken. Conjugacy problem is 
known to be hard in some groups and hence it seems to be a
reasonable assumption at this moment. There is another worry: if Alice sends
$k$ and $\delta$ separately then one can find $a$ from the equation $x=(\delta
k)(a\gamma)^{-1}$, since $\gamma$ is computable. However, this is circumvented
easily by sending the product $\delta k$ not $\delta$ and $k$ individually and
keeping $k$ random. 
\section{An interesting family of $p$-groups}
It is well known that cyclic groups have abelian automorphism groups. The first
person to give an example of a non-abelian group with an abelian automorphism
groups is G.A. Miller in \cite{miller} which was generalized by Struik in
\cite{struik}. There are three non-abelian groups with abelian automorphism
group in the Hall-Senior table \cite{halltable}, they are 
nos.~91, 92 and 99. Miller's example is no.~99. In \cite{jamali}, Jamali
generalized nos.~91 and 92. His generalization of no.~91 is in one direction,
it increases the exponent of the group. 

Jamali in the same paper generalizes group no.~92 in two directions, the
size of the exponent and the number of generators. 
His generalization was restrictive
in that it works only for the prime
$2$. There are other examples of families of Miller $p$-groups in the literature,
the most notable one is the family of $p$-groups, for an arbitrary prime $p$,
given by Jonah and 
Konisver in \cite{jonah}. This was generalized to an arbitrary number of
generators by Earnley in \cite{earnley}. There are other examples by Martha
Morigi in \cite{morigi} and Heineken and Liebeck in \cite{hl}. All these
examples of Miller groups given in \cite{earnley,hl,jonah,morigi} are special
groups, i.e., the commutator and the center are the 
same. For special groups the key exchange protocols do not
work as noted earlier. So there is no Miller $p$-group, readily
available in the literature, for arbitrary prime $p$ which can be 
used right away in construction of the protocol. The only other source are
groups nos.~ 91, 92 and 99 in the Hall Senior table \cite{halltable} and their
generalizations, notice that these groups are not special but are
$2$-groups. Of the three generalizations, the generalization of no.~ 92 best
fits our
criterion because it is generalized in two directions, \textit{viz}.~
number of generators and exponent of the center and moreover it is not
special; $Z(G)=A\times G^\prime$ where $A$ is a cyclic group. So once we
generalize it for arbitrary primes, it has ``three degrees of freedom'', the
number of generators, exponent of center and the prime; which makes it
attractive for cryptographic purposes. 
 
In the rest of the section we use Jamali's definition in \cite{jamali} to
define a family  
of $p$-groups for arbitrary prime. So this family is a generalization of
Jamali's example and assuming transitivity of generalizations, ultimately a
generalization of group no.~ 92 in the Hall-Senior table \cite{halltable}. We
study automorphisms of this 
group and show that the group is Miller if and only if $p=2$, but this family
of groups always have an abelian central automorphism group which is fairly
large. We then attempt to 
build a key exchange protocol as described earlier using the central
automorphisms. We start with the definition of the group $G_n(m,p)$.
\begin{definition}
Let $G_n(m,p)$ be a group generated by $n+1$ elements 

\noindent $\{a_0,a_1,a_2,\ldots ,a_n\}$ where $p$ is a prime number
and $m\geq 2$ and $n\geq 3$ are integers. The group is
defined by the following relations: 
\[a_1^p=1,\;\;\;a_2^{p^m}=1,\;\;\;a_i^{p^2}=1\;\;\; \text{for}\;\;\; 
3\leq i\leq n,\;\;\;a_{n-1}^p=a_0^p.\]
\[[a_1,a_0]=1,\;\;\;[a_n,a_0]=a_1,\;\;\;[a_{i-1},a_0]=a_i^p
\;\;\text{for}\;\;\;3\leq i\leq n.\]
\[[a_i,a_j]=1 \;\;\;\text{for}\;\;\;1\leq i<j\leq n.\] 
\end{definition}
We state some facts about the group $G_n(m,p)$ whose proof is by direct
computation (see \cite[Section 2.9]{ayan1}).
\begin{description}
\item[a] $G_n(m,p)^\prime$ the derived subgroup of $G_n(m,p)$ is an
  elementary abelian group
 $\langle a_1,a_3^p,\ldots a_n^p\rangle\simeq\mathbb{Z}_p^{n-1}$. 
\item[b] $Z(G_n(m,p))=\langle a_2^p\rangle \times G^\prime$.
\item[c] $G_n(m,p)$ is a $p$-group of class $2$.
\item[d] $G_n(m,p)$ is a PN group.
\end{description}
\begin{proposition}
$G_n(m,p)$ is a polycyclic group and every element of $g\in G_n(m,p)$ can be
uniquely expressed in the form
$g=a_0^{\alpha_0}a_1^{\alpha_1}a_2^{\alpha_2}a_3^{\alpha_3}\ldots
a_n^{\alpha_n}$, where\\
$0\leq\alpha_i <p$ for $i=0,1$;  $0\leq\alpha_2 <p^m$, 
$0\leq\alpha_i <p^2$ for $i=3,4,\ldots ,n$.
\end{proposition}
\begin{proof}
Define $G_0=G_n(m,p)=\langle a_0,a_1,a_2,\ldots ,a_n\rangle$,
$G_1=\langle a_1,a_2,\ldots a_n\rangle$ and similarly $G_k=\langle
a_k,a_{k+1},\ldots ,a_n\rangle$ for $k\leq n$. Since $G_1$ is a finitely
generated 
abelian group, it is a polycyclic group \cite[Proposition 3.2]{sims}. It
is fairly straightforward to show that 
\[G_1\triangleright
G_2\triangleright \ldots\triangleright G_n\triangleright\langle 1\rangle\]
is a polycyclic series and $\{a_1,\ldots ,a_n\}$ a polycyclic generating
sequence of $G_1$.

It is easy to see from the relations of the group that $G_1$ is normal in
$G_0$ and $G_0/G_1$ 
is cyclic. It is also straightforward to show that $\langle
a_iG_{i+1}\rangle=G_i/G_{i+1}$ and $|a_iG_{i+1}|=|a_i|$ and hence any element
of the group has a unique representation of the above form. We would call an
element represented in the above form a \emph{collected word}. See also \cite[Chapter
9, Proposition 4.1]{sims}.    
\end{proof}

\noindent\textbf{Computation with $G_n(m,p)$}: Our group $G_n(m,p)$ 
is of class 2, i.e., commutators of weight 3 are identity, computations
become real nice and easy. Let us demonstrate the product of two collected
words  
$g=a_0^{\alpha_0}a_1^{\alpha_1}a_2^{\alpha_2}a_3^{\alpha_3}a_4^{\alpha_4}$ and 
$h=a_0^{\beta_0}a_1^{\beta_1}a_2^{\beta_2}a_3^{\beta_3}a_4^{\beta_4}$.
To compute $gh$ we use concatenation and form the word
$a_0^{\alpha_0}a_1^{\alpha_1}a_2^{\alpha_2}a_3^{\alpha_3}a_4^{\alpha_4}a_0^{\beta_0}a_1^{\beta_1}a_2^{\beta_2}a_3^{\beta_3}a_4^{\beta_4}$  
and note that $a_i$'s commute except for $a_0$ hence one tries to move $a_0$
towards the left using the identity 
\[a_ia_0=a_0a_i[a_i,a_0]=\left\{
\begin{array}{ccc}
a_0a_ia_{i+1}^p &\text{for}& 1\leq i<n\\
a_0a_ia_1 &\text{for}& i=n.
\end{array}
\right.
\]  
Further note, since commutators are in the center of the group,
$a_{i+1}^p$ or $a_1$ can be moved anywhere. Once $a_0$ is moved to the extreme
left 
the word formed is the collected word of $gh$. This process is often
referred to in the literature as \emph{collection}. Computing the inverse of an
element can be similarly done.

We now prove that the group of central automorphisms of the group
$G_n(m,p)$ for an arbitrary prime $p$ is abelian. For sake of simplicity we
denote $G_n(m,p)$ by $G$ for the rest of the article, and use notations from
Theorem \ref{adney}. 

\begin{lemma}
In $G$, $\mathcal{R}=Z(G)=\mathcal{K}$.
\end{lemma}
\begin{proof}
Using the notation from Theorem \ref{adney}, we see that in $G$, $a=m-1$,
$b=1$ and $c=m$ hence $d=m-1$. Clearly, $\mathcal{R}=Z(G)$ hence
$\mathcal{K}\subseteq Z(G)$. 

\noindent Let $x\in Z(G)$, if $x\in G^\prime$ then height$(xG^\prime)=\infty$
and we 
are done. If not, then $x=z_1z_2$ where $z_1\in\langle a_2^p\rangle$ and
$z_2\in G^\prime$. Then $xG^\prime=z_1G^\prime$ and hence
height$(xG^\prime)\geq 1$. 
\end{proof}
It is easy to see that $\mathcal{R}/G^\prime=Z(G)/G^\prime=\langle
a_2^pG^\prime\rangle$ and hence from Theorem \ref{adney} we prove the
following theorem:
\begin{theorem}
$\text{Aut}_c(G)$ is abelian.
\end{theorem}
\subsection{Automorphisms of $G_n(m,p)$}In this section we describe the
automorphisms of groups of this kind. The discussion is, in more than one way,
an adaptation of the work of Jamali \cite{jamali} and generalizes his
main theorem. 
\begin{lemma}
Let $x=a_0^{\beta_0} a_1^{\beta_1}a_2^{\beta_2}\ldots a_n^{\beta_n}$, where
$\beta_i$, $i=0,1,2\ldots ,n$ are integers be an element of $G$. If $p=2$ then
${\beta_0}$ is $1$ 
and 
\begin{itemize}
\item[\textbullet] $x^2=a_1^{\beta_n} a_2^{2\beta_2}a_3^{\gamma_3}\ldots
a_{n-2}^{\gamma_{n-2}}a_{n-1}^{\gamma_{n-1}+2}a_n^{\gamma_n}$ for $p=2$. Where
$\gamma_i=2(\beta_{i-1}+\beta_i)$. 
\item[\textbullet] $x^p=a_2^{p\beta_2}a_3^{p\beta_3}\ldots
  a_{n-2}^{p\beta_{n-2}}a_{n-1}^{p\beta_{n-1}+p{\beta_0}}a_n^{p\beta_n}$ for
  $p\neq 2$.
\end{itemize}
\end{lemma}
\begin{proof}
For the case $p=2$ we just collect terms and use the relation
$a_{n-1}^2=a_0^2$.

\noindent For $p\neq 2$ using Proposition 8.3 we have 
\begin{eqnarray*}
\lefteqn{x^p=(a_0^{\beta_0} a_1^{\beta_1}a_2^{\beta_2}\ldots
a_{n-1}^{\beta_{n-1}}a_n^{\beta_n})^p}\\
&=&(a_0^{\beta_0})^p(a_1^{\beta_1}a_2^{\beta_2}\ldots 
a_{n-1}^{\beta_{n-1}}a_n^{\beta_n})^p\\
&=&a_0^{p{\beta_0}}a_2^{p\beta_2}a_3^{p\beta_3}\ldots a_n^{p\beta_n}
\end{eqnarray*}
Using the relation $a_{n-1}^p=a_0^p$ we have  
$$a_0^{p{\beta_0}}a_2^{p\beta_2}a_3^{p\beta_3}\ldots
a_n^{p\beta_n}=a_2^{p\beta_2}a_3^{p\beta_3}\ldots a_{n-2}^{p\beta_{n-2}}a_{n-1}^{p\beta_{n-1}+p{\beta_0}}a_n^{p\beta_n}$$  
\end{proof}
\noindent For the group $G$ we note that $H=\langle a_1,a_2,a_3,\ldots
a_n\rangle $ is the maximal 
abelian normal subgroup of $G$ and is characteristic. It follows that
the $H^p$ is also characteristic. Following \cite{jamali}, 
we define two decreasing sequences of characteristic subgroups
$\{K_i\}_{i=0}^{n-1}$ such that $$K_0=H \;\text{and}\;
K_i/K_{i-1}^p=Z(G/K_{i-1}^p)\;\;(1\leq i\leq n-1)$$
and
$\{L_i\}$ such that
$$ L_0=H\;\text{and}\;L_i=\{h:\;h\in H,\;h^p\in[G,L_{i-1}]\}\;(1\leq i\leq
n-1)$$ 
It follows easily that 

$$K_i=\langle a_1,a_2,\ldots ,a_{n-i},a^p_{n-i+1},\ldots ,a_n^p\rangle \;
1\leq i\leq n-1$$ 
$$L_1=\langle a_1,v,a_3,\ldots ,a_n\rangle $$
$$L_i=\langle a_1,v,a_3^p,\ldots ,a^p_{i+1},a_{i+2},\ldots
,a_n\rangle\;\;\;2\leq i\leq n-1$$ 
where $v=a_2^{p^{m-1}}$. For $3\leq i\leq n$ we have
$$K_{n-i}\cap L_{i-2}=\langle a_1,v,a_3^p,\ldots
a_{i-1}^p,a_i,a_{i+1}^p,\ldots a_n^p\rangle =\langle v,a_i,G^\prime\rangle .$$ 
Also $K_{n-2}\cap L_0=\langle a_2,G^\prime\rangle $. 

Since $\langle v,a_i,G^\prime\rangle$ and $\langle a_2,G^\prime\rangle$ are
characteristic, for any $\theta\in\text{Aut}(G)$,
\begin{eqnarray*}
\theta(a_2)&=a_2^{k_2}z&\text{where}\;\;z\in
G^\prime\;\;\text{and}\;\;k_2\in\mathbb{N}\\
\theta(a_i)&=a_i^{k_i}v^{r_i}z&\text{where}\;\;z\in G^\prime;\;
k_i\in\mathbb{N};\;\;i=3,4,\ldots,n;\;\;0\leq r_i<p. 
\end{eqnarray*}
It is clear that not all $k_2$ and $k_i$ will make $\theta$ an
automorphism. To begin with, if $\theta$ is an automorphism then
$\gcd(k_i,p)=1$ for 
all $k_i$, and we may choose $k_i$, such that $0<k_i<p$ for $i=3,4,\ldots,n$.

Let $\theta(a_0)=a_0^{\beta_0}a_1^{\beta_1}a_2^{\beta_2}\ldots
a_n^{\beta_n}$. Since $\theta(a_0^p)=\theta(a_{n-1}^p)=\theta
(a_{n-1})^p=a^{pk_{n-1}}_{n-1}$, from Lemma 10.4
$$a_{n-1}^{pk_{n-1}}=a_2^{p\beta_2}a_3^{p\beta_3}\ldots
a_{n-2}^{p\beta_{n-2}}a_{n-1}^{p\beta_{n-1}+p{\beta_0}}a_n^{p\beta_n}\;\;\text{for}\;\;p\neq
2$$ implying $\beta_0+\beta_{n-1}\equiv k_{n-1}\mod p$, $p^{m-1}|\beta_2$
and $p|\beta_i$ for $i=3,4,\ldots,n-2,n$. Hence
$\theta(a_0)=a_0^{k_0}a_{n-1}^{\beta_{n-1}}v^rz$ where $0\leq r<p$. We
changed $\beta_0$ to $k_0$ to maintain
uniformity in notations. 

Notice the relation $[a_i,a_0]=a_{i+1}^p$ for $i=2,3,\ldots,n$ implies that 
\[[\theta(a_i),\theta(a_0)]=\theta(a_{i+1})^p=a_{i+1}^{pk_{i+1}}.\]
It follows that
$[a_i^{k_i},a_0^{k_0}a_{n-1}^{\beta_{n-1}}]=a_{i+1}^{pk_{i+1}}$ which is the
same as $[a_i^{k_i},a_0^{k_0}]=a_{i+1}^{pk_{i+1}}$, which implies that
$[a_i,a_0]^{k_0k_i}=a_{i+1}^{pk_{i+1}}$. Recall that $G$ is a 
$p$-group of class $2$ and $a_{n-1}$ commutes with $a_i$ for $i\geq 2$. From
these we have a recursive formula for $k_i$, (also see \cite[Theorem 2.9.7]{ayan1}):
choose $k_0$ such that $0<k_0<p$ and $k_2$ such that $0<k_2<p^m$ and
$\gcd(k_2,p)=1$ and then define $k_{i+1}=k_0k_{i}\mod p$ for
$i=2,3,4,\ldots,(n-1)$; and $k_1=k_0k_n\mod p$.  In \cite[Proposition 2.3]{jamali} Jamali proves that
for $p=2$, all automorphisms of $G$ are central. We have just proved that for
$p\neq 2$ there is a non central automorphism, take $k_0>1$; the
following theorem follows from Theorem 8.4.
\begin{theorem}
The group $G_n(m,p)$ is Miller if and only if $p=2$.
\end{theorem} 
\subsection{Description of the Central Automorphisms}
Notice that $G$ is a PN group, so there is a one-one correspondence between
$\text{Aut}_c(G)$ and $\text{Hom}(G,Z(G))$. Since it is known from our
earlier discussion that $Z(G)=\langle
a_2^p\rangle\times G^\prime$, $\text{Hom}(G,Z(G))=\text{Hom}(G,\langle
a_2^p\rangle)\times\text{Hom}(G,G^\prime)$. It follows:
$\text{Aut}_c(G)=A\times B$ where
\[A=\{\sigma\in\text{Aut}_c(G):x^{-1}\sigma(x)\in\langle a_2^p\rangle\}\]
\[B=\{\sigma\in\text{Aut}_c(G):x^{-1}\sigma(x)\in G^\prime\}\] 
Elements of $A$ can be explained in a very nice way. Pick
a random integer $k$ such that $k=lp+1$ where $0\leq l<p^{m-1}$ and a
random 
subset $R$ (could be empty) of $\{0,3,4,\ldots n\}$, and then an arbitrary
automorphism in $A$ is 
\begin{eqnarray}\label{oneeqn}
\nonumber &\sigma(a_1)=a_1\\
\nonumber&\sigma(a_2)=a_2^k\\
\nonumber&\sigma(a_i)=\left\{
\begin{array}{cccc}
\nonumber a_i &\text{if}\;\;i \not\in R\\
\nonumber a_i\left(a_2^{p^{m-1}}\right)^{r_i}&\text{if}\;\; i\in R 
\end{array}\right.\\
\end{eqnarray}
We use indexing in $\{0,3,4,\ldots ,n\}$ to order $R$ and $0<r_i<p$ is an
integer corresponding to $i\in R$.
Conversely, any element in $A$ can be described this way. It follows
from the definition of $A$ that $$|A|=p^{m-1}\times p^{n-1}=p^{m+n-2}.$$

\noindent The automorphism $\phi\in B$ is of the form
\begin{eqnarray}\label{twoeqn}
\phi(x)=\left\{
\begin{array}{llll}
a_1&\text{if}& x=a_1\\
a_iz&\text{if}& x=a_i,\;\;i\in\{0,2,3,\ldots ,n\}
\end{array}\right.
\end{eqnarray}
where $z\in G^\prime$. 

We note that $\dfrac{G}{Z(G)}$ is an abelian group and hence
$\text{Inn}(G)$ is abelian and hence
$\text{Inn}(G)\subseteq\text{Aut}_c(G)$. We further note from the
commutator relations in $G$ that
$\text{Inn}(G)\subseteq B$.
\subsection{Using these automorphisms in key-exchange protocol I}
Let us briefly recall the key-exchange protocol described before. Alice and Bob
decide on a group $G$ and a non-central element $g\in G\setminus Z(G)$ in
public. Alice then chooses an arbitrary automorphism $\phi_A$ and sends Bob
$\phi_A(g)$. Similarly Bob picks an arbitrary automorphism $\phi_B$ and sends
Alice $\phi_B(g)$. Since the automorphisms commute, both of them can compute
$\phi_A(\phi_B(g))$, which is their private key. The most devastating attack on
the system is the one in which Oscar, looking at $g$, $\phi_A(g)$ and
$\phi_B(g)$, can predict what $\phi_A(\phi_B(g))$
will look like, i.e., a GDHP attack.
\begin{definition}[Parity condition for elements in $G$]
If $g=a_0^{\beta_0} a_1^{\beta_1}a_2^{\beta_2}a_3^{\beta_3}\ldots
a_n^{\beta_n}$ is an arbitrary element of $G$, i.e., $0\leq{\beta_0} <p$,
$0\leq 
\beta_1<p$, $0\leq\beta_2 <p^m$ and $0\leq \beta_i <p^2$ for $3\leq i\leq
n$. Then the vector $v:=({\beta_0},\beta_3,\beta_4,\ldots ,\beta_n)$ is called
the parity 
of $g$. Two elements $g$ and $g^\prime$ are said to be of the same parity condition
if $v=v^\prime\mod p$, where $v^\prime$ is the parity of $g^\prime$.
\end{definition}
\begin{lemma}
If $\phi:G\rightarrow G$ is any central automorphism then $g$ and
$\phi(g)$ have the same parity condition for any $g\in G$. 
\end{lemma}
\begin{proof}
Notice that an automorphism $\phi$ either belongs to $A$ or $B$ or is of the
form $\phi(g)=gf_{\phi}(g)g_{\phi}(g)$ where $f_{\phi}\in\text{Hom}(G,Z(G))$
and $g_{\phi}\in\text{Hom}(G,G^\prime)$. So we might safely ignore
elements from $A$, since they only affect the exponent of $a_2$. Also note that
$a_1$ being in the commutator remains fixed under any central automorphism.

\noindent So we need to be concerned with elements of $B$, from the
description of $B$, and each commutator is a word in $p$-powers of the
generators and from the fact that $G^\prime\subset Z(G)$, the lemma follows. 
\end{proof}
Now let us understand what an element in $A$ does to an element
$g\in G$. We use notations from Equation \ref{oneeqn}.
\begin{lemma}
Let
$g=a_0^{\beta_0}a_1^{\beta_1}a_2^{\beta_2}a_3^{\beta_3}\ldots 
a_n^{\beta_n}$,  $\phi\in A$ and if

$\phi(g)=a_0^{\beta_0^\prime}a_1^{\beta_1^\prime}a_2^{\beta_2^\prime}a_3^{\beta_3^\prime}\ldots  
a_n^{\beta_n^\prime}$ then $\beta_i=\beta_i^\prime$ for $i\neq 2$ and

$\beta_2^\prime=k\beta_2+p^{m-1}\sum\limits_{i\in R}r_i\beta_i \mod p^m$.
\end{lemma}
\begin{proof}
Notice that from Equation \ref{oneeqn}, it is clear that elements of $A$ only
affect the exponent of $a_2$, so $\beta_i^\prime =\beta_i$ for $i\neq 2$
follows trivially. From the definition of $A$ and simple computation it follows
that $\beta_2^\prime=k\beta_2+p^{m-1}\sum\limits_{i\in R}r_i\beta_i\mod p^m$. 
\end{proof} 
In the key exchange protocol I, we will only use automorphisms
from\footnote{In light of Lemma 10.6, we believe that adding automorphisms
  from $B$ is not going to add to the security of the system.} $A$. As noted earlier there are two kinds of attack, GDLP (the discrete
logarithm 
problem in automorphisms) and GDHP (the Diffie-Hellman problem in
automorphisms). We have earlier stated that GDLP is equivalent to
finding 
the automorphism from the action of the automorphism on one element. It seems
that for one to find the automorphism discussed in the previous lemma, one has
to find $k$, $R$ and $r_i$. Notice that
$\beta_2^\prime=k\beta_2+p^{m-1}\sum\limits_{i\in R}r_i\beta_i\mod p^m$, is a
knapsack in $\beta_2$ and $p^{m-1}$. Solving that knapsack is not enough to
compute the image of any element, because $R$ is not known so $\beta_i$'s are
not known. We shall show in a moment that the security
of the key exchange protocol depends on the difficulty of this knapsack,
 but solving this knapsack does not help Oscar to find the
 automorphism, just partial information about the automorphism comes out. 

Next we show that though it seems to be secure under GDLP, but if the knapsack
is solved then the system is broken by GDHP. This proves that GDHP is a weaker
problem than GDLP in $G_n(m,p)$.
 Let
 $g=a_0^{\beta_0}a_1^{\beta_1}a_2^{\beta_2}a_3^{\beta_3}\hspace*{\fill}\\
\ldots a_n^{\beta_n}$, 
then as discussed before for $\phi,\psi\in\text{Aut}_c(G)$, with notation
from Equation 2 and $k_i\in\mathbb{N}$ for $i=3,4,\ldots,n$:
\begin{equation}\label{eqn1}
\phi(g)=a_0^{\beta_0}a_1^{\beta_1}a_2^{k_2\beta_2+p^{m-1}\sum\limits_{i\in
  R}r_i\beta_i}a_3^{\beta_3+k_3p}\ldots a_n^{\beta_4+k_4p}
\end{equation}
\begin{equation}\label{eqn2}
\psi(g)=a_0^{\beta_0}a_1^{\beta_1}a_2^{k_2^\prime\beta_2+p^{m-1}\sum\limits_{i\in R^\prime}r_i^\prime\beta_i}a_3^{\beta_3+k_3^\prime p}\ldots a_n^{\beta_4+k_4^\prime
  p} 
\end{equation}
From direct computation it follows
that the exponent of $a_2$ in
$\phi(\psi(g))$ is 
\begin{equation}\label{eqn3}
k_2\left(k_2^\prime\beta_2+p^{m-1}\sum\limits_{i\in 
  R^\prime}r_i^\prime\beta_i\right)+p^{m-1}\sum\limits_{i\in
  R}r_i\beta_i
\end{equation} where $k_2=lp+1$ and $k_2^\prime=l^\prime p+1$, $0\leq
  l,l^\prime <p^{m-1}$.
The exponent of $a_0, a_1$ stays the same and the exponent of $a_i$ will be
$\beta_i+(k_i+k_i^\prime)p\mod p^2$ for $3\leq i\leq n$. As mentioned before
since we are using only automorphisms from $A$, i.e., $\phi$ and $\psi$ are in
$A$ hence $k_i=k_i^\prime=0$ for $i=3,4,\ldots,n$.

Notice that $g$, Equations \ref{eqn1} and \ref{eqn2} are public, so Oscar sees
those. Since the
exponents of $a_0,a_1,a_3,\ldots,a_n$ are predictable, the key Alice
and Bob want to establish is the exponent of $a_2$ in 
$\phi\left(\psi(g)\right)$, which is given by Equation 6. Since Oscar sees
Equations \ref{eqn1} and \ref{eqn2}, if he can compute $k_2$ from
$k_2\beta_2+p^{m-1}\sum\limits_{i\in R}r_i\beta_i\mod p^m$, then he can
compute $p^{m-1}\sum\limits_{i\in R}r_i\beta_i$ and the scheme is
broken. But, $k_2=lp+1$ for some $l\in[0,p^{m-1})$ hence
$$k_2\beta_2+p^{m-1}\sum\limits_{i\in R}r_i\beta_i\mod p^m$$ reduces to
$$\beta_2+lp\beta_2+p^{m-1}\sum\limits_{i\in R}r_i\beta_i\mod p^m.$$ Since
$\beta_2$ is public, Oscar can compute
$lp\beta_2+p^{m-1}\sum\limits_{i\in R}r_i\beta_i\mod p^m$. Notice that finding
$k_2$ is equivalent to finding $l$, hence one of the security assumptions is
that there is no polynomial time algorithm to find $l$ from
\begin{equation}lp\beta_2+p^{m-1}\sum\limits_{i\in R}r_i\beta_i\mod
  p^m.
\end{equation} 
Let us write 
\begin{equation}\label{eqn4}
M=lp\beta_2+p^{m-1}\sum\limits_{i\in R}r_i\beta_i\mod p^m, 
\end{equation}
then
\begin{equation*}
M=lp\beta_2\mod p^{m-1}.
\end{equation*} 
Now, if $lp<p^{m-1}$ and $\gcd(\beta_2,p)=1$, then one can find $lp$
from the above equation and the scheme is broken. So the only hope of
making a secure cryptosystem out of key exchange protocol I and
the group $G_n(m,p)$ is to take $l=kp^{m-2}$ where
$k=0,1,2,\ldots,(p-1)$. In this case, if we set $l=lp^{m-2}$ and
$l^\prime=l^\prime p^{m-2}$ in Equation \ref{eqn3}, then the key will be
\begin{eqnarray*}
\lefteqn{\left(1+lp^{m-1}\right)\left((1+l^\prime
  p^{m-1})\beta_2+p^{m-1}\sum\limits_{i\in
  R^\prime}r_i^\prime\beta_i\right)+p^{m-1}\sum\limits_{i\in
  R}r_i\beta_i}\\
&=\left(1+lp^{m-1}+l^\prime p^{m-1}\right)\beta_2+p^{m-1}\sum\limits_{i\in R^\prime}r_i^\prime\beta_i+p^{m-1}\sum\limits_{i\in
  R}r_i\beta_i\mod p^m\\
&=\left(\left(1+lp^{m-1}\right)\beta_2+p^{m-1}\sum\limits_{i\in
  R}r_i\beta_i\right)+l^\prime p^{m-1}\beta_2+p^{m-1}\sum\limits_{i\in
  R^\prime}r_i^\prime\beta_i\mod p^m
\end{eqnarray*}

Now the information in the last equation is easy to compute from the
public information, Equations \ref{eqn1} and \ref{eqn2}; so the Key
Exchange Protocol I is broken for automorphisms from $A$ of $G_n(m,p)$
when $\gcd(\beta_2,p)=1$. 

Now if $\gcd(p,\beta_2)\neq 1$, i.e., $\beta_2=kp^i$ for some
$i\in[1,m)$ and $1\leq k <p$, then an attack similar to the above breaks the
system. The insight behind these attacks is that any solution to Equation
\ref{eqn4} can be thought of as the image of $g$ under an automorphism
$\phi^\prime\in A$. We are talking about a solution to Equation
\ref{eqn4}, which is easy to find, for which $\phi^\prime(g)=M$
and then Shpilrain's attack breaks the system. 
\section{Implementation}
There is not much reason left to go into the details of an implementation.
We briefly mention that this cryptosystem can be implemented without
any reference to the group $G_n(m,p)$. Once the element
$g=a_0^{\beta_0}a_1^{\beta_1}a_2^{\beta_2}\ldots a_n^{\beta_n}$ is
fixed, Alice can send Bob $k_2\beta_2+p^{m-1}\sum\limits_{i\in R}r_i\beta_i
\mod p^m$ and similarly Bob can send Alice
$k_2^\prime\beta_2+p^{m-1}\sum\limits_{i\in R^\prime}r_i^\prime\beta_i  \mod
p^m$. Since Alice and Bob know their
own $k_2$, $\sum\limits_{i\in
  R}r_i\beta_i$ and $k_2^\prime$, $\sum\limits_{i\in
  R^\prime}r_i^\prime\beta_i$ respectively, they can both compute the
private key or the shared secret. Since the only operation involved in
computing 
the private key is multiplication and addition$\mod p^m$, there
can be a very fast implementation of this cryptosystem.
\section{Conclusion}
In this paper we studied a key exchange protocol using commuting
automorphisms in a non-abelian $p$-group. Since any nilpotent group is a direct
product of its Sylow subgroups, the study of
nilpotent groups can be reduced 
to the study of $p$-groups. We argued that our study
is a generalization of the Diffie-Hellman 
key exchange and is a generalization of the discrete log problem. Other
public key systems like the El-Gamal cryptosystem which uses the
discrete logarithm problem
is adaptable to our methods. This is the first attempt to
generalize the discrete logarithm problem in the way we did. 

We should try to find other groups and try our system in terms of GDLP and
GDHP. As we noted earlier, GDHP is a subproblem of the GDLP, and we saw in
$G_n(m,p)$, GDHP is a much easier problem than GDLP. 
Our example was of the form $d>b$ in Theorem \ref{adney}. The next step is to
look at groups where $d=b$. We note from Theorem \ref{sanders}, if a $p$-group
$G$ is a PN group then $\text{Aut}_c(G)$ is a $p$-group and since $p$-groups
have nontrivial centers; one can work in that center with our scheme. In this
case we would be generalizing to arbitrary nilpotentcy class while still working
with central automorphisms. 

Lastly we note that, if we were using some representation for this finitely
presented group $G$, for example, matrix representation of the group over a
finite field $\mathbb{F}_q$; then security of the system in $G_n(m,p)$ becomes
the discrete logarithm problem in a matrix algebra
\cite{menezes1,menezes2}. Since the discrete 
logarithm problem in matrices is only as secure as the discrete 
logarithm problem in finite fields, there is no known advantage to go
for matrix representation, but there might 
be other representations of interest. 

There is one conjecture that comes out of this work and
we end with that.   
\begin{conjecture}
If $G$ is a Miller $p$-group for an odd prime $p$, then $G$ is special.
\end{conjecture}
\noindent\textbf{Acknowledgment}:
The author wishes to thank Fred Richman and Rustam Stolkin; they read
the whole manuscript and made 
valuable suggestions. The author is indebted to the referee for his
kind comments.
\nocite{*}
\bibliography{Implementation}

\begin{thebibliography}{10}

\bibitem{adney}
A.~Adney and T.~Yen.
\newblock Automorphisms of p-group.
\newblock {\em Illinois Journal of Mathematics}, 9:137--143, 1965.

\bibitem{braid3}
I.~Anshel, M.~Anshel, B.~Fisher, and D.~Goldfield.
\newblock New key agreement protocols in braid group cryptography.
\newblock In {\em CT-RSA 2001}, number 2020 in Lecture Notes in Computer
  Science, pages 1--15. Springer, 2001.

\bibitem{braid4}
I.~Anshel, M.~Anshel, and D.~Goldfeld.
\newblock An algebraic method for public-key cryptography.
\newblock {\em Math. Research Letters}, 6:287--291, 1999.

\bibitem{links}
J.~S. Birman.
\newblock {\em Braids, Links and Mapping Glass Groups}.
\newblock Number~82 in Annals of Mathematics Studies. Princeton University
  Press, 1974.

\bibitem{blakeCurve}
I.~Blake, G.~Seroussi, and N.~Smart.
\newblock {\em Eliptic Curves in Cryptography}.
\newblock Number 265 in London Mathematical Society, Lecture Note Series.
  Cambridge University Press, 1999.

\bibitem{blake1}
I.~F. Blake and T.~Garefalakis.
\newblock On the complexity of the {D}iscrete {L}ogarithm and
  {D}iffie-{H}ellman problems.
\newblock {\em Journal of Complexity}, 20, 2004.

\bibitem{boneh}
D.~Boneh.
\newblock The {D}ecision {D}iffie-{H}ellman problem.
\newblock In {\em Algorithmic number theory (Portland, OR, 1998)}, number 1423
  in Lecture Notes in Computer Science, pages 48--63. Springer, Berlin, 1998.

\bibitem{polybraid}
J.~H. Cheon and B.~Jun.
\newblock An polynomial time algoritm for tha braid {D}iffie-{H}ellman
  conjugacy problem.
\newblock In {\em Advances in cryptography -- CRYPTO 2003}, number 2729 in
  Lecture Notes in Computer Science, pages 212--225. Springer, Berlin, 2003.

\bibitem{curran}
M.~Curran.
\newblock Semidirect product groups with abelian automorphism groups.
\newblock {\em J. Austral. Math, Soc.}, Series A(42):84--91, 1987.

\bibitem{braid5}
P.~Dehornoy.
\newblock Braid-based cryptogrpahy.
\newblock {\em Contemporary Mathematics}, 360:1--33, 2004.

\bibitem{earnley}
B.~E. Earnley.
\newblock {\em On finite Groups whose group of automorphisms is abelian}.
\newblock PhD thesis, Wayne State University, 1975.

\bibitem{elgamal}
T.~Elgamal.
\newblock A public key cryptosystem and a signature scheme based on discrete
  logarithms.
\newblock {\em Lecture Notes in Comput. Sci., 196, Springer, Berlin.}, pages
  10--18, 1985.

\bibitem{fournelle}
T.~Fournelle.
\newblock Elementary abelian $p$-groups as automorphism group of infinite
  group.
\newblock {\em I. Math. Z.}, 167:259--270, 1979.

\bibitem{galb1}
S.~Galbraith and V.~Rotger.
\newblock Easy decision {D}iffie-{H}ellman groups.
\newblock {\em LMS Journal of Computation and Mathematics}, 7:201--218, 2004.

\bibitem{GAP4}
The GAP~Group.
\newblock {\em {GAP -- Groups, Algorithms, and Programming, Version 4.3}},
  2002.
\newblock \verb+(http://www.gap-system.org)+.

\bibitem{halltable}
M.~Hall and J.~Senior.
\newblock {\em The groups of order $2^n\;\;(n\leq 6)$}.
\newblock Macmillan, 1964.

\bibitem{hall}
P.~Hall.
\newblock {\em The Edmonton notes on nilpotent groups}.
\newblock Queen Mary college mathematics notes, Cambridge, 1969.

\bibitem{hl}
H.~Heineken and H.~Liebeck.
\newblock The occurrence of finite groups in the automorphism group of
  nilpotent groups of class $2$.
\newblock {\em Archives of Mathematics}, 25:8--16, 1974.

\bibitem{hopkins}
C.~Hopkins.
\newblock Non-abelian groups whose groups of isomorphism are abelian.
\newblock {\em Ann. of Math}, 29(1-4):508--520, 1927.

\bibitem{jamali}
A.-R. Jamali.
\newblock Some new non-abelian 2-groups with abelian automorphism groups.
\newblock {\em Journal of Group Theory}, 5:53--57, 2002.

\bibitem{jonah}
D.~Jonah and M.~Konvisser.
\newblock Some non-abelian $p$-groups with abelian automorphism groups.
\newblock {\em Archives of Mathematics}, 26:131--133, 1975.

\bibitem{kaplansky}
I.~Kaplansky.
\newblock {\em Infinite Abelian Groups}.
\newblock The University of Michigan Press, 1969.

\bibitem{khukhro}
E.~Khukhro.
\newblock {\em $p$-Automorphisms of finite $p$-groups}.
\newblock Number 246 in London Mathematical Society, Lecture Note Series.
  Cambridge University Press, 1997.

\bibitem{braid2}
K.~H. Ko, D.~H. Choi, M.~S. Cho, and J.~W. Lee.
\newblock New signature scheme using conjugacy problem.
\newblock http://eprint.iacr.org/2002/168, 2002.

\bibitem{braid1}
K.~H. Ko, S.~J. Lee, J.~H. Cheon, J.~W. Han, J.~sung Kang, and C.~Park.
\newblock New public-key cryptosystem using braid groups.
\newblock In M.~Bellare, editor, {\em Advances in Cryptology -- CRYPTO 2000},
  number 1880 in Lecture Notes in Computer Science, pages 166--183, 2000.

\bibitem{koblitz1}
N.~Koblitz.
\newblock {\em A course in number theory and cryptography}.
\newblock Number 114 in Graduate Texts in Mathematics. Springer-Verlag, New
  York, second edition, 1994.

\bibitem{koblitz}
N.~Koblitz.
\newblock {\em Algebraic aspects of cryptography}.
\newblock Number~3 in Algorithms and Computation in Mathematics.
  Springer-Verlag, Berlin, 1998.

\bibitem{parameters}
N.~Koblitz, A.~Menezes, and S.~Vanstone.
\newblock The state of {E}lliptic {C}urve {C}ryptography.
\newblock {\em Designs, Codes and Cryptography}, 19:173--193, 2000.

\bibitem{kurosh}
A.~Kurosh.
\newblock {\em The Theory of Groups}, volume 1 \& 2.
\newblock Chelsea Publishing Company, 1960.

\bibitem{ayan1}
A.~Mahalanobis.
\newblock {\em Diffie-{H}ellman Key Exchange Protocol, Its Generalization and
  Nilpotent Groups}.
\newblock PhD thesis, Florida Atlantic University, August 2005.
\newblock http://eprint.iacr.org/2005/223.

\bibitem{ayanJP}
A.~Mahalanobis.
\newblock Abelian groups, homomorphisms and central automorphisms of nilpotent
  groups.
\newblock {\em JP Journal of Algebra, Number Theory and Applications},
  7(1):69--81, 2007.

\bibitem{menezes1}
A.~J. Menezes and S.~A. Vanstone.
\newblock A note on cyclic group, finite fields and discrete logarithm problem.
\newblock {\em Applicable Algebra in Engineering, Communication and Computing},
  3(1):67--74, 1992.

\bibitem{menezes2}
A.~J. Menezes and Y.-H. Wu.
\newblock The discrete logarithm problem in {${\rm GL}(n,q)$}.
\newblock {\em Ars Combinatoria}, 47:23--32, 1997.

\bibitem{miller}
G.~Miller.
\newblock A non-abelian group whose group of isomorphism is abelian.
\newblock {\em Messenger Math.}, 43:124--125, 1913.

\bibitem{morigi}
M.~Morigi.
\newblock On $p$-groups with abelian automorphism group.
\newblock {\em The Mathematical Journal of the University of Padova},
  92:47--58, 1994.

\bibitem{rotman}
J.~J. Rotman.
\newblock {\em An introduction to the theory of groups}.
\newblock Springer-Verlag, 1994.

\bibitem{sanders}
P.~R. Sanders.
\newblock The central automorphism of a finite group.
\newblock {\em J. London Math. Soc.}, 44:225--228, 1969.

\bibitem{scott}
W.~Scott.
\newblock {\em Group Theory}.
\newblock Dover, 1964.

\bibitem{igor}
I.~E. Shparlinski.
\newblock Security of polynomial transformations of {D}iffie-{H}ellman key.
\newblock {\em Finite fields and their applications}, 10:123--131, 2004.

\bibitem{sh1}
V.~Shpilrain and G.~Zapata.
\newblock Combinatorial group theory and public key cryptography.
\newblock {\em Applicable Algebra in Engineering, Communication and Computing},
  17(3-4):291--302, 2006.

\bibitem{sims}
C.~Sims.
\newblock {\em Computation with finitely presented groups}.
\newblock Cambridge University Press, Cambridge, 1994.

\bibitem{stinson}
D.~Stinson.
\newblock {\em Cryptography: Theory and Practice}.
\newblock CRC Press, 2 edition, 2002.

\bibitem{struik}
R.~R. Struik.
\newblock Some non-abelian 2-groups with abelian automorphism groups.
\newblock {\em Archives of Mathematics}, 39:299--302, 1982.

\bibitem{gonzalez}
M.~I.~G. Vasco and I.~E. Sharlinski.
\newblock On the security of {D}iffie-{H}ellman bits.
\newblock In {\em Cryptography and computational number theory}, Progress in
  Computer Science and Applied Logic, pages 257--268. Birkh\"{a}user, Basel,
  2001.

\bibitem{maria}
M.~S. Voloshina.
\newblock {\em On the holomorph of a discrete group}.
\newblock PhD thesis, University of Rochester, 2003.

\bibitem{zeanhaus}
H.~Zassenhaus.
\newblock {\em The theory of groups}.
\newblock Chelsea, New York, 1958.

\end{thebibliography}
\bibliographystyle{abbrv}
\end{document}